\documentclass[12pt,reqno]{amsart}
\usepackage{amssymb,latexsym,amsmath,mathrsfs, xypic} 
\usepackage[pdftex, colorlinks=false, pdfborder={1 1 1}, pdfborderstyle={/S/U/W 1}]{hyperref} 

\newcommand{\K}{\mathcal{K}}
\renewcommand{\L}{\mathcal{L}}
\newcommand{\U}{\mathcal{U}}
\newcommand{\M}{{\mathcal{M}}}
\newcommand{\C}{\mathbb{C}}
\newcommand{\N}{\mathbb{N}}

\newcommand{\Z}{\mathbb{Z}}
 
\newcommand{\R}{\mathbb{R}}

\newcommand{\Od}{{\mathcal{O}_2}}
\newcommand{\Oinf}{{\mathcal{O}_\infty}}
\newcommand{\Td}{{\mathcal{T}_2}}

\newcommand{\ot}{\otimes}
\newcommand{\pf}{\noindent {\mbox{\textit{Proof}. }} }
\newcommand{\ie}{\textit{i.e.}\,\,} 
\newcommand{\eg}{\textit{e.g.\/}\ }
\newcommand{\cst}{C$^*$}

\newtheorem{thm}{Theorem}[section]
\newtheorem{cor}[thm]{Corollary}

\newtheorem{prop}[thm]{Proposition}
\newtheorem{defi}[thm]{Definition}

\theoremstyle{remark}
\newtheorem{rem}[thm]{Remark}
\newtheorem{rems}[thm]{Remarks}

\newtheorem{ques}[thm]{Question}

\numberwithin{equation}{section}
\title{Continuous fields of properly infinite \cst-algebras
}

\author{Etienne Blanchard}



\begin{document} 
\begin{abstract} 
Any unital separable continuous $C(X)$-algebra with properly infinite fibres is properly infinite 
as soon as the compact Hausdorff space $X$ has finite topological dimension. 
We study conditions under which this is still the case 
in the infinite dimensional case. 
\end{abstract} 

\maketitle 

\section{Introduction} 
One of the basic \cst-algebras studied in the classification programme 
launched by G. Elliott (\cite{Ell94}) 
of nuclear \cst-algebras 
through K-theoretical invariants 
is the Cuntz \cst-algebra $\Oinf$ generated by infinitely many isometries 
with pairwise orthogonal ranges (\cite{Cun77}). 
This \cst-algebra is pretty rigid in so far as it is a strongly self-absorbing \cst-algebra (\cite{TW07}): 
Any separable unital continuous $C(X)$-algebra $A$ 
the fibres of which are isomorphic to the same strongly self-absorbing \cst-algebra $D$ 
is a trivial $C(X)$-algebra 
provided the compact Hausdorff base space $X$ has finite topological dimension. 
Indeed, the strongly self-absorbing \cst-algebra $D$ tensorially absorbs 
the JiangÐ-Su algebra~$\mathcal{Z}$ (\cite{Win09}). 
Hence, this \cst-algebra $D$ is K$_1$-injective (\cite{Ror04}) and 
the $C(X)$-algebra $A$ satisfies $A\cong D\ot C(X)$ (\cite{DW08}). 
\\ \indent 
But I. Hirshberg, M. R\o rdam and W. Winter have built 
a non-trivial unital continuous \cst-bundle 
over the infinite dimensional compact product $\Pi_{n=0}^\infty\, S^2$ 
such that all its fibres are isomorphic to 
the strongly self-absorbing UHF algebra of type $2^\infty$ 
(\cite[Example 4.7]{HRW07}). 
More recently, M. D\u{a}d\u{a}rlat has constructed in \cite[\S 3]{Dad09} 
for all pair $(\Gamma_0,\Gamma_1)$ of countable abelian torsion groups 
a unital separable continuous $C(X)$-algebra $A$ such that 
\begin{enumerate} 
\item[--] the base space $X$ is the compact Hilbert cube $X=B_\infty$ of infinite dimension, 
\item[--] all the fibres $A_x$ ($x\in B_\infty$) are isomorphic to 
the strongly self-absorbing Cuntz \cst-algebra $\Od$ 
generated by two isometries $s_1, s_2$ satisfying $1_\Od=s_1^{}s_1^*+s_2^{}s_2^*\,$, 
\item[--] $K_i(A)\cong C(Y_0, \Gamma_i)$ for $i=0, 1\,$, 
where $Y_0\subset [0,1]$ is the canonical Cantor set. 
\end{enumerate} 
These $K$-theoretical conditions imply that the $C(B_\infty)$-algebra $A$ is not a trivial one.
But these arguments do not work anymore 
when the strongly self-absorbing algebra $D$ is the Cuntz algebra $\Oinf$ (\cite{Cun77}), 
in so far as $K_0(\Oinf)=\mathbb{Z}$ is a torsion free group. 

We describe in this article the link between several notions of proper infiniteness for $C(X)$-algebras 
which appeared during the recent years (\cite{KR00},\cite{BRR08}, \cite{CEI08}, \cite{RR11}).  
We then study whether  
certain unital continuous $C(X)$-algebras with fibres $\Oinf$, 
especially the Pimsner-Toeplitz algebra (\cite{Pim95}) of Hilbert $C(X)$-modules 
with fibres of dimension greater than $2$ 
and with compact base space $X$ of infinite topological dimension, 
are properly infinite \cst-algebras. 

\medskip 
I especially thank E. Kirchberg and a referee for a few inspiring remarks. 

\section{A few notations} 

We present in this section the main notations which are used in this article. 
We denote by $\N=\{0, 1, 2, \ldots\}$ the set of positive integers and 
we denote by $[S]$ the closed linear span of a subset $S$ in a Banach space. 

\begin{defi}(\cite{Dix69}, \cite{Kas88}, \cite{Blan97}) Let $X$ be a compact Hausdorff space and 
let $C(X)$ be the \cst-algebra of continuous function on $X\,$. 
\begin{enumerate} 
\item[--] A unital $C(X)$-algebra is a unital \cst-algebra $A$ 
endowed with a unital morphism of \cst-algebra from $C(X)$ to the centre of $A$. 
\item[--] For all closed subset $F\subset X$ and all element $a\in A$, 
one denotes by $a_{|F}$ the image of $a$ in the quotient $A_{|F}:=A/C_0(X\setminus F)\cdot A$. 
If $x\in X$ is a point in $X$, one calls fibre at $x$ the quotient $A_x:=A_{|\{ x\}}$ and 
one write $a_x$ for $a_{|\{ x\}}$. 
\item[--] The $C(X)$-algebra $A$ is said to be continuous if 
the upper semicontinuous map $x\in X\mapsto\| a_x\|\in\mathbb{R}_+$ is continuous 
for all $a\in A$. 
\end{enumerate}
\end{defi}

\begin{rems} 
a) (\cite{Cun81}, \cite{BRR08}) For all integer $n\geq 2$, 
the \cst-algebra $\mathcal{T}_n:=\mathcal{T}(\C^n)$ is the universal unital \cst-algebra 
generated by $n$ isometries $s_1, \ldots, s_n$ satisfying the relation 
\begin{equation}\label{Tn} 
s_1^{}s_1^* +\ldots+s_n^{}s_n^*\leq 1\;.
\end{equation} 
\noindent b) A unital \cst-algebra $A$ is said to be \textit{properly infinite} if and only if 
one the following equivalent conditions holds true (\cite{Cun77}, \cite[Proposition 2.1]{Ror03}): 
\begin{enumerate} 
\item[--] the \cst-algebra $A$ contains two isometries with mutually orthogonal range projections, 
\ie $A$ unitally contains a copy of $\Td\,$, 
\item[--] the \cst-algebra $A$ contains a unital copy of the simple Cuntz \cst-algebra $\Oinf$ 
generated by infinitely many isometries with pairwise orthogonal ranges. 
\end{enumerate} 
c) If $A$ is a \cst-algebra and $E$ is a Hilbert $A$-module, 
one denotes by $\L(E)$ the set of adjointable $A$-linear operators acting on $E$ (\cite{Kas88}). 
\end{rems} 

\section{Global proper infiniteness} 
Proposition 2.5 of \cite{BRR08} and section 6 of \cite{Blan13} entail 
the following stable proper infiniteness 
for continuous $C(X)$-algebras with properly infinite fibres. 
 \begin{prop} 
 Let $X$ be a second countable perfect compact Hausdorff space, \ie without any isolated point, and 
 let $A$ be a separable unital continuous $C(X)$-algebra with properly infinite fibres. 

\noindent 1) There exist: 
\begin{enumerate} 
\item[(a)] a finite integer $n\geq 1\,$, 
\item[(b)] a covering $X= \mathop{F_1}\limits^o\cup\ldots\cup\mathop{F_n}\limits^o$ 
by the interiors of closed balls $F_1, \ldots, F_n\,$, 
\item[(c)] unital embeddings of \cst-algebra $\sigma_k: \Oinf\hookrightarrow A_{| F_k}$ 
($1\leq k\leq n\,$). 
\end{enumerate} 
2) The tensor product $M_p(\mathbb{C})\otimes A$ is properly infinite 
for all large enough integer $p$. 
\end{prop} 
\pf 
1) For all point $x\in X$, 
the semiprojectivity of the \cst-subalgebra $\Oinf\hookrightarrow A_x$ 
(\cite[Theorem 3.2]{Blac04}) entails that 
there are a closed neighbourhood $F\subset X$ of the point $x$ and 
a unital embedding $\Oinf\ot C(F)\hookrightarrow A_{| F}$ of $C(F)$-algebra. 
The compactness of the topological space $X$ enables to conclude. 

\medskip\noindent 
2) Proposition 2.7 of \cite{BRR08} entails that the \cst-algebra $M_{2^{n-1}}(A)$ is properly infinite. 
Proposition~2.1 of \cite{Ror97} implies that $M_p(A)$ for all integer $p\geq 2^{n-1}$.
\qed

\medskip\begin{rem} 
If $X$ is a second countable compact Hausdorff space and 
$A$ is a separable unital continuous $C(X)$-algebra, 
then $\widetilde{X}:=X\times [0,1]$ is a perfect compact space, 
$\widetilde{A}:=A\otimes C([0, 1])$ is a unital continuous $C(\widetilde{X})$-algebra and 
every morphism of unital \cst-algebra $\Oinf\to\widetilde{A}$ induces 
a unital $\ast$-homomorphism $\Oinf\to A$ 
by composition with the projection map $\widetilde{A}\to A$ 
coming from the injection $\quad x\in X\mapsto (x,0)\in\widetilde{X}\,$. 
\end{rem} 

\medskip 
The proper infiniteness of the tensor product $M_p(\mathbb{C})\otimes A$ 
does not always imply that the \cst-algebra $A$ is properly infinite (\cite{HR98}). 
Indeed, there exists a unital \cst-algebra $A$ which is not properly infinite, 
but such that the tensor product $M_2(\mathbb{C})\otimes A$ is a properly infinite \cst-algebra 
(\cite[Proposition 4.5]{Ror03}). 
The following corollary nevertheless holds true. 

\begin{cor}\label{cor3.3} 
Let $\jmath_0, \jmath_1$ denote the two canonical unital embeddings 
of the Cuntz extension $\Td$ in the full unital free product $\Td\ast_\C\Td$ and 
let $\tilde{u}\in\U(\Td\ast_\C\Td)$ be a K$_1$-trivial unitary satisfying 
$\jmath_1(s_1)=\tilde{u}\cdot\jmath_0(s_1)$ (\cite[Lemma 2.4]{BRR08}). 

The following assertions are equivalent: 
\begin{enumerate}
\item[(a)] The full unital free product $\Td\ast_\C\Td$ is $K_1$-injective. 
\item[(b)] The unitary $\tilde{u}$ belongs to 
the connected component $\U_0(\Td\ast_\C\Td)$ of $1_{\Td\ast_\C\Td}$. 
\item[(c)] Every separable unital continuous $C(X)$-algebra $A$ with properly infinite fibres 
is a properly infinite \cst-algebra. 
\end{enumerate} 
\end{cor} 
\pf (a)$\Rightarrow$(b) 
A unital \cst-algebra $A$ is called K$_1$-injective if and only if 
all $K_1$-trivial unitaries $v\in\U(A)$ are homotopic to the unit $1_A$ in $\U(A)$ 
(see \eg \cite{Roh09}). 
Thus, (b) is a special case of (a) 
since $K_1(\Td\ast_\C\Td)=\{1\}$ (see \eg \cite[Lemma 4.4]{Blan10}). 

\smallskip\noindent (b)$\Rightarrow$(c) 
Let $A$ be a separable unital continuous $C(X)$-algebra with properly infinite fibres. 
Take a finite covering $X= \mathop{F_1}\limits^o\cup\ldots\cup\mathop{F_n}\limits^o$ 
such that there exist unital embeddings $\sigma_k:\Td\to A_{|F_k}$ for all $1\leq k\leq n$. 
Set $G_k:=F_1\cup\ldots\cup F_k\subset X$ 
and let us construct by induction isometries $w_k\in A_{|G_k}$ such that 
the two projections $w_kw_k^*$ and $1_{|G_k}-w_kw_k^*$ are 
properly infinite and full in the restriction $A_{|G_k}\,$: 

\noindent 
-- If $k=1$, the isometry $w_1:=\sigma_1(s_1)$ 
has the requested properties. 

\noindent 
-- If $k\in\{1, \ldots, n-1\}$ and the isometry $w_k\in A_{|G_k}$ is already constructed, 
then Lemma~2.4 of \cite{BRR08} implies that there exists 
a morphism of unital \cst-algebra $\pi_k:\Td\ast_\C\Td\to A_{|G_k\cap F_{k+1}}$ 
satisfying 
\begin{equation}\label{3.1} 
\begin{array}{l}
-\quad \pi_k(\jmath_0(s_1))=w_k{}_{|G_k\cap F_{k+1}}\,,\\
-\quad \pi_k(\jmath_1(s_1))=\sigma_{k+1}(s_1){}_{|G_k\cap F_{k+1}}=
\pi_k(\tilde{u})\cdot w_k{}_{|G_k\cap F_{k+1}}\;. 
\end{array} 
\end{equation} 
If the unitary $\tilde{u}$ belongs to the connected component $\U_0(\Td\ast_\C\Td)$, 
then $\pi_k(\tilde{u})$ is homotopic 
to $1_{A_{|G_k\cap F_{k+1}}}=\pi_k(1_{\Td\ast_\C\Td})$ 
in $\U(A_{|G_k\cap F_{k+1}})$, 
so that $\pi_k(\tilde{u})$ admits a unitary lifting $z_{k+1}$ in $\U_0(A_{|F_{k+1}})$ 
(see \eg \cite[ Lemma 2.1.7]{LLR00}). 
The only isometry $w_{k+1}\in A_{|G_{k+1}}$ satisfying the two constraints 
\begin{equation}\label{form4.2 } 
\begin{array}{l}
-\quad w_{k+1}{}_{| G_k}=w_k\,,\\
-\quad w_{k+1}{}_{| F_{k+1}}=(z_{k+1})^*\cdot\sigma_{k+1}(s_1) 
\end{array} 
\end{equation} 
verifies that the two projections $w_{k+1}w_{k+1}^*$ and $1_{|G_{k+1}}-w_{k+1}w_{k+1}^*$ are 
properly infinite and full in $A_{|G_{k+1}}\,$. 

The proper infiniteness of the projection $w_n^{}w_n^*$ in $A_{|G_n}=A$ implies that 
the unit $1_A=w_n^* w_n^{}=w_n^*\cdot w_nw_n^*\cdot w_n$ is also 
a properly infinite projection in $A$, 
\ie the \cst-algebra $A$ is properly infinite. 

\smallskip\noindent (c)$\Rightarrow$(a) 
The \cst-algebra $\mathcal{D}\!:=\!\!\{f\in C([0, 1] , \Td\ast_\C\Td)\,;\, 
f(0)\in\jmath_0(\Td)\,\mathrm{and}\,f(1)\!\in\!\jmath_1(\Td)\,\}$ is 
a unital continuous $C([0, 1])$-algebra the fibres of which are all properly infinite. 
Thus, condition (c) implies that the \cst-algebra $\mathcal{D}$ is properly infinite, 
a statement which is equivalent to the K$_1$-injectivity of $\Td\ast_\C\Td$ 
(\cite[Proposition 4.2]{Blan10}). \qed 

\section{A question of proper infiniteness} 
We describe in this section the link between the different notions of proper infiniteness 
which have been introduced during the last decades. 

\medskip 
Recall first the link between properly infinite \cst-algebras and 
properly infinite Hilbert \cst-modules studied by K. T. Coward, G. Elliott and C. Ivanescu. 
\begin{prop}(\cite{CEI08}) Let $A$ be a stable \cst-algebra and let $a\in A_+$ be a positive element. 
The following assertions are equivalent:
\begin{enumerate}
\item[(a)] $a$ is properly infinite in $A$, \ie $a\oplus a\precsim a$ in $\K\ot A$ 
(\cite[definition 3.2]{KR00}). 
\item[(b)] There is an embedding of Hilbert $A$-module 
$\ell^2(\N)\ot A\hookrightarrow [aA]$ (\cite{CEI08}). 
\end{enumerate} 
\end{prop} 
\pf (a) $\Rightarrow$ (b) If there exists an infinite sequence $\{ d_i\}$ in $A$ such that 
$a=d_i^*d_i^{}\geq\sum_{j\in\N} d_j^{}d_j^*$ for all $i\in\N$, 
then $[aA]\supset\sum_j[d_jA]\cong\ell^2(\N)\ot A$. 

\smallskip\noindent (b) $\Rightarrow$(a) 
One has embeddings of Hilbert $A$-modules $[aA]\oplus [aA]\subset \ell_2(\N)\ot A\subset [aA]$. 
\qed 

\medskip 
The following holds for continuous deformations of properly infinite \cst-algebras. 
\begin{prop}\label{prop4-2} 
Let $X$ be a second countable compact Hausdorff space, 
let $H$ be a separable Hilbert $C(X)$-module with non-zero fibres and 
let $a\in\K(H)$ be a strictly positive compact contraction. 
Consider the following assertions: 
\begin{enumerate} 
\item[(a)] All the Hilbert spaces $H_x$ are properly infinite, \ie 
the operator $a_x$ is properly infinite in the \cst-algebra $\K(H)_x\cong\K(H_x)$ for all $x\in X$. 
\item[(b)] $H$ is a properly infinite Hilbert $C(X)$-module, \ie $a$ is properly infinite in $\K(H)$. 
\item[(c)] The multiplier \cst-algebra $\L(H)=\M(\K(H) )$ is properly infinite. 
\end{enumerate} 
Then $(c)\Rightarrow(b)\Rightarrow(a)$. 
But $(a)\not\Rightarrow(b)$ and $(b)\not\Rightarrow(c)$.
\end{prop} 
\pf (c)$\Rightarrow$(b) 
If $\sigma:\Td=C^*<s_1, s_2>\,\to \L(H)$ is a unital $\ast$-homomorphism, then 
the two elements $d_1=\sigma(s_1)\cdot a^{1/2}$ and $d_2=\sigma(s_2)\cdot a^{1/2}$ 
satisfy $d_i^*d_j^{}=\delta_{i=j}\cdot a$ in $\K(H)$. 

\medskip\noindent(b)$\Rightarrow$(a) 
The relations $c_i^*c_j^{}=\delta_{i=j}\cdot a$ between $3$ operators $c_1, c_2, a$ 
in a $C(X)$-algebra $D$ entails that 
$(c_i)_x^*(c_j)_x^{}=\delta_{i=j}\cdot a_x$ in the quotient $D_x=D/ C_0(X\setminus\{ x\})\cdot D$ 
for all $x\in X$. 

\medskip\noindent(a)$\not\Rightarrow$(b) 
Let $B_3=\{ (x_1, x_2, x_3)\in\R^3\,;\,x_1^2+x_2^2+x_3^2\leq 1\}$ be the unit ball of dimension~$3$, 
let $B_3^+, B_3^-$ be the two open semi-disks 
$B_3^+=\{(x_1, x_2, x_3)\in B_3 \,;\, x_3>-\frac{1}{2}\}$, 
$B_3^-=\{(x_1, x_2, x_3)\in B_3\,;\, x_3<\frac{1}{2}\}$ 
and 
let $S_2=\{ (x_1, x_2, x_3)\in B_3\,;\,\,x_1^2+x_2^2+x_3^2=1\}\subset B_3$ be the unit sphere of dimension $2$. 
\\ \indent 
The self-adjoint operator $f\in C(B_3)\ot M_2(\C)\cong C(B_3, M_2(\C) )$ given by 
\begin{equation}
f(x_1, x_2, x_3)=\frac{1}{2}\cdot\left(
\begin{array}{cc}1+x_3&x_1-\imath\, x_2\\x_1+\imath\, x_2&1-x_3\end{array}\right) 
\end{equation} 
is a positive contraction since each self-adjoint matrix $f(x_1, x_2, x_3)\in M_2(\C)$ satisfies 
\begin{equation} 
f(x_1, x_2, x_3)^2=f(x_1, x_2, x_3)+\frac{1}{4}(x_1^2+x_2^2+x_3^2-1)\cdot    1_{M_2(\C)}\;,
\end{equation} 
\ie $(f(x_1, x_2, x_3)-\frac{1}{2}\cdot  1_{M_2(\C)})^2=
\frac{x_1^2+x_2^2+x_3^2}{4}\cdot    1_{M_2(\C)}\leq(\frac{1}{2}\cdot    1_{M_2(\C)})^2\;$. 
\\ \indent 
The nontrivial Hilbert $C(B_3)$-module 
$F:= [f\cdot \left(\begin{smallmatrix}C(B_3)\\ C(B_3)\end{smallmatrix}\right)]$ 
satisfies the two isomorphisms of Hilbert $C(B_3)$-module: 
\begin{equation} 
\begin{array}{l}
-\quad F\cdot C_0(B_3^+)\cong C_0(B_3^+)\oplus C_0(B_3^+\setminus S_2\cap B_3^+)\\ 
-\quad F\cdot C_0(B_3^-)\cong C_0(B_3^-)\oplus C_0(B_3^-\setminus S_2\cap B_3^-)\;.  
\end{array}
\end{equation} 

\smallskip 
The set $B_\infty:=\{ x\in\ell^2(\N)\,;\,\sum_p |x_p|^2\leq 1\}$ 
is a metric compact space called the \textit{complex Hilbert cube} 
when equipped with the distance $d((x_p), (y_p))=\sum_p\, 2^{-p-2}\, |x_p-y_p|$. 
Denote by $E_{DD}$ the non-trivial Hilbert $C(B_\infty)$-module with fibres $\ell^2(\N)$ 
constructed by J. Dixmier and A. Douady (\cite[\S 17]{DD63}, \cite[Proposition 3.6]{BK04a}). 

\smallskip 
Finally, consider the product $X := B_\infty\times B_3$ and the Hilbert $C(X)$-module 
\begin{equation} 
H:=\,E_{DD}\otimes C(B_3)\; \oplus\; C(B_\infty)\otimes F \;.
\end{equation} 
The two Hilbert $C(X)$-submodules $H\cdot C_0(B_\infty\times B_3^+)$ and 
$H\cdot C_0(B_\infty\times B_3^-)$ are properly infinite, 
\ie there exist embeddings of Hilbert $C(X)$-module 
$\ell^2(\N)\ot C_0(B_\infty\times B_3^+)\hookrightarrow H\cdot C_0(B_\infty\times B_3^+)$ and 
$\ell^2(\N)\ot C_0(B_\infty\times B_3^-)\hookrightarrow H\cdot C_0(B_\infty\times B_3^-)$. 
Hence, all the fibres of the Hilbert $C(X)$-module $H$ are properly infinite Hilbert spaces, 
\ie $\ell^2(\N)\hookrightarrow H_x$ for all point $x$ in the compact space 
$X=B_\infty\times B_3^+\cup B_\infty\times B_3^-$ (\cite{CEI08}). 
But the Hilbert $C(X)$-module $H$ is not properly infinite 
\ie $\ell^2(\N)\ot C(X)\not\hookrightarrow H$ (\cite[Example 9.11]{RR11}). 
The equality $C(B_3^{})=C_0(B_3^+)+C_0(B_3^-)$ only implies that 
$\ell^2(\N)\ot C(X)\hookrightarrow H\oplus H$. 

\medskip\noindent(b)$\not\Rightarrow$(c) 
There exists a continuous field $\tilde{H}$ of Hilbert spaces 
over the compact space Y:=$B_\infty\times (B_3)^\infty$ such that 
$\tilde{H}=[a\cdot \tilde{H}]$ for some properly infinite contraction $a\in\K(\tilde{H})$ and 
the \cst-algebra $\L(\tilde{H})$ is not properly infinite.  
Indeed, let $\eta\in\ell^\infty(B_\infty, \ell^2(\N)\oplus\C)$ be the section 
$x\mapsto x\oplus \sqrt{1-\| x\|^2}$, 
let $\ddot{F}$ be the closed Hilbert $C(B_\infty)$-module 
$\ddot{F}:=[C(B_\infty, \ell^2(\N)\oplus 0)+C(B_\infty)\cdot\eta]$, 
let $\theta_{\eta, \eta}\in\L(\ddot{F})$ be the projection $\zeta\mapsto\eta\langle\eta,\zeta\rangle$ and 
let $E_{DD}=(1-\theta_{\eta, \eta})\cdot \ddot{F}$ be the Hilbert $C(B_\infty)$-submodule 
built in \cite{DD63}. 
Define also the sequence of contractions $\widetilde{f}=(\widetilde{f}_n)$ 
in $\ell^\infty\Bigl(\N, M_2\bigl(C((B_3)^\infty)\bigr)\Bigr)$ by 
\begin{equation}
x_n=(x_{n, k})\in(B_3)^\infty\longmapsto \widetilde{f}_n(x_n):=f(x_{n, n})\in M_2(\C)\,. 
\end{equation} 
The Hilbert $C(Y)$-module 
\begin{equation} 
\tilde{H}:=C(Y)\,\oplus\, E_{DD}\ot C((B_3)^\infty)\,\oplus\, 
C(B_\infty)\ot[\widetilde{f}\cdot
\ell^2\Bigl(\N, \left(\begin{smallmatrix}C( (B_3)^\infty)\\ C( (B_3)^\infty)\end{smallmatrix}\right)\,\Bigr)] 
\end{equation} 
has the desired properties (\cite[Example 9.13]{RR11}). 
\qed 
\begin{ques} 
-- What can be written when the operator $a$ in Proposition~\ref{prop4-2} is 
a projection in a unital continuous $C(X)$-algebra? \\
-- Is the full unital free product $\Td\ast_\C\Td$ a properly infinite \cst-algebra 
which is not $K_1$-injective? 
(see the equivalence (a)$\Leftrightarrow$(c) in Corollary~\ref{cor3.3}) 
\end{ques} 
\section{The Pimsner-Toeplitz algebra of a Hilbert $C(X)$-module} 
We look in this section at the proper infiniteness of 
the unital continuous $C(X)$-algebras with fibres $\Oinf$ 
corresponding to the Pimsner-Toeplitz $C(X)$-algebras 
of Hilbert $C(X)$-modules with infinite dimension fibres. 

\medskip 
\begin{defi}(\cite{Pim95}) 
Let $X$ be a compact Hausdorff space and 
let $E$ be a full Hilbert $C(X)$-module, \ie without any zero fibre. \\ 
a) The full Fock Hilbert $C(X)$-module $\mathcal{F}(E)$ of $E$ is 
the direct sum 
\begin{equation}\label{2.1}
\mathcal{F}(E):=\mathop{\bigoplus}\limits_{m\in\N}\;E^{(\otimes_{C(X)})\,m}\;,
\end{equation} 
where $E^{(\otimes_{C(X)})\,m}:=\left\{
\begin{array}{ll}
C(X)&\mathrm{if}\; m=0\,, \\
E\otimes_{C(X)}\ldots\otimes_{C(X)} E\;\; (m\,\mathrm{terms}) &\mathrm{if}\; m\geq 1\,.
\end{array}\right.$ 
\\ b) The Pimsner-Toeplitz $C(X)$-algebra $\mathcal{T}(E)\,$ of $E$
is the unital subalgebra of 
the $C(X)$-algebra $\L(\mathcal{F}(E)\,)$ 
of adjointable $C(X)$-linear operators acting on $\mathcal{F}(E)$ 
generated by the creation operators $\ell(\zeta)$ ($\zeta\in E$), 
where 
\begin{equation}\label{2.2}
\begin{array}{llcll}
-\quad \ell(\zeta)\,(f \cdot\hat{1}_{C(X)})&\!\!:=f\cdot\zeta=\zeta\cdot f &
\textrm{for}&f\in C(X)&\mathrm{ and}\\ 
-\quad \ell(\zeta)\,(\zeta_1\otimes\ldots\otimes\zeta_k)&
\!\!:=\zeta\otimes\zeta_1\otimes\ldots\otimes\zeta_k &
\textrm{for}&\zeta_1,\ldots, \zeta_k\in E& \mathrm{if}\,k\geq 1\,.
\end{array}
\end{equation}
c) Let $(C^*(\mathbb{Z}), \Delta)$ be the abelian compact quantum group 
generated by a unitary $\mathbf{u}$ with spectrum the unit circle and 
with coproduct $\Delta(\mathbf{u})=\mathbf{u}\otimes\mathbf{u}$. 
Then, there is a unique coaction $\alpha$ 
of the Hopf \cst-algebra $(C^*(\mathbb{Z}), \Delta)$ 
on the Pimsner-Toeplitz $C(X)$-algebra $\mathcal{T}(E)$ 
such that $\alpha\bigl(\ell(\zeta)\bigr)=\ell(\zeta)\otimes\mathbf{u}$ 
for all $\zeta\in E$, \ie 
\begin{equation} 
\begin{array}{cccccc}\alpha:&
\mathcal{T}(E)&\to&\mathcal{T}(E)\otimes C^*(\mathbb{Z})&=&C(\mathbb{T}, \mathcal{T}(E))\\
&\ell(\zeta)&\mapsto&\ell(\zeta)\otimes \mathbf{u}\quad&=&(z\mapsto\ell(z\zeta)\,)
\end{array} 
\end{equation} 
\end{defi} 

The fixed point $C(X)$-subalgebra 
$\mathcal{T}(E)^{\alpha}=\{ a\in\mathcal{T}(E)\,;\,\alpha(a)=a\otimes 1\}$ 
under this coaction is the closed linear span 
\begin{equation}\label{Cor3.2}
\mathcal{T}(E)^{\alpha}=
\Bigl[C(X).1_{\mathcal{T}(E)}+\sum_{k\geq 1}\,\ell( E)^k\cdot \left(\ell( E)^k\right)^*\Bigr]\,.
\end{equation}
\indent 
Besides, the projection $P\in\L(\,\mathcal{F}(E) )$ onto the submodule $E$ 
induces a quotient morphism of $C(X)$-algebra 
$a\in\mathcal{T}(E)^\alpha\mapsto 
\mathbf{\overline{q}}(a):=P\cdot a\cdot P\in \mathcal{K}(E)+C({X})\cdot 1_{\L(E)}\subset\L(E)$. 

\medskip 
\begin{prop}\label{prop4.2} 
Let $X$ be a second countable compact Hausdorff perfect space 
and let $E$ be a separable Hilbert $C(X)$-module with infinite dimensional fibres. 

\noindent 1) There exist a covering 
$X=\mathop{F_1}\limits^o\cup\ldots\cup\mathop{F_m}\limits^o$ 
by the interiors of closed subsets $F_1, \ldots, F_m$ 
and $m$ norm $1$ sections $\zeta_1, \ldots, \zeta_m$ in $E$ such that 
$\mathcal{T}(E)={C^*\!\!<\!\mathcal{T}(E)^\alpha, \ell(\zeta_1), \ldots, \ell(\zeta_m)\!>}$, 
$\left(\ell(\zeta_k)^{}\ell(\zeta_k)^*\right)_{|F_k}\,$ and 
$\left(1-\ell(\zeta_k)^{}\ell(\zeta_k)^*\right)_{|F_k}$ 
are properly infinite projections in $\mathcal{T}(E)_{|F_k}$ 
for all index $k\in\{ 1,\ldots, m\}$\,. 

\noindent 2) Set $G_k:=F_1\cup\ldots\cup F_k$ for all integer $k\in\{1, \ldots, m\}$ and 
$\bar{G}_l:=G_l\cap F_{l+1}$ for all integer $l\in\{1, \ldots, m-1\}$. 
If $\xi(l)\in E_{|G_l}$ is a section such that $\|\xi(l)_y\|=1$ for all $y\in\bar{G}_l$, 
then there is a unitary $w_l\in\mathcal{T}(E)_{|\bar{G}_l}$ such that 
\begin{enumerate} 
\item[(a)] $w_l\cdot \ell(\xi(l))_{|\bar{G}_l}=\ell(\zeta_{l+1})_{|\bar{G}_l}\,$, 
\item[(b)] $w_l\oplus 1_{|\bar{G}_l}$ is homotopic to $1_{|\bar{G}_l}\oplus 1_{|\bar{G}_l}$ 
among the unitaries in $M_2(\mathcal{T}(E)_{|\bar{G}_l})\,$. 
\end{enumerate} 

\noindent 3) If for all K$_1$-trivial unitary $w_k\in\mathcal{T}(E)_{|\bar{G}_k}$ 
there is a unitary $z_{k+1}\in\mathcal{T}(E)^\alpha{}_{| F_{k+1}}$ 
such that $(z_{k+1})_{| \bar{G}_k}=w_k$ ($1\leq k\leq m-1$), 
then there is a section $\xi\in E$ satisfying 
\begin{equation}
\forall\, x\in X\,,\quad \|\xi_x\|=1\;,
\end{equation} 
so that Lemma 6.1 of \cite{Blan13} implies that the \cst-algebra $\mathcal{T}(E)$ is properly infinite. 
\end{prop} 
\pf 
1) Given a point $x\in X$ and a unit vector $\zeta\in E_x$, 
let $\xi_1, \xi_2, \xi_3$ be three norm $1$ sections in $E$ such that \
$(\xi_1)_x=\zeta$ and 
the matrix $a=[\langle\xi_i, \xi_j\rangle]\in M_3( C(X) )$ satisfies $a_x=1_3\in M_3(\C)$. 
Let $F\subset X$ be a closed neighbourhood of $x$ such that $\| a_y-1_3\|\leq 1/2$ for all $y\in F$. 
Define the sections $\xi_1', \xi_2', \xi_3'$ in $E_{| F}$ by 
\begin{equation}\label{5.6} 
\left(\begin{smallmatrix}(\xi_1')_y\\(\xi_2')_y\\(\xi_3')_y\end{smallmatrix}\right)=
\left(\begin{smallmatrix}(\xi_1)_y\\(\xi_2)_y\\(\xi_3)_y\end{smallmatrix}\right) \cdot 
(a_y^*\, a_y^{})^{-1/2} \quad \quad\mathrm{for}\,\mathrm{ all}\, y\in F\,. 
\end{equation} 
One has $\Bigl[\langle\xi_i', \xi_j'\rangle\Bigr]=
(a_{|F}^*\, a_{|F}^{})^{-1/2}\cdot (a_{|F}^*\, a_{|F}^{})\cdot (a_{|F}^*\, a_{|F}^{})^{-1/2}=1$ in $M_3( C(F) )$. 
Hence, $\ell(\xi_1')^{}\ell(\xi_1')^*$ and $q:=1_{|F}-\ell(\xi_1')^{}\ell(\xi_1')^*$ are 
properly infinite projections in $\mathcal{T}(E)_{| F}$ since 
\begin{equation}\label{5,7}
\begin{array}{l}
-\quad 1_{|F}-q=\ell(\xi_1')^{}\ell(\xi_1')^*\geq
\ell(\xi_1')^{}\ell(\xi_2')^{}\ell(\xi_2')^*\ell(\xi_1')^*+\ell(\xi_1')^{}\ell(\xi_3')^{}\ell(\xi_3')^*\ell(\xi_1')^*\\ 
-\quad q\geq \ell(\xi_2')^{}\ell(\xi_2')^*+\ell(\xi_3')^{}\ell(\xi_3')^*\geq 
\ell(\xi_2')^{}\,q^2\,\ell(\xi_2')^*+\ell(\xi_3')^{}\,q^2\,\ell(\xi_3')^*\,,
\end{array}
\end{equation} 
so that there exist unital $\ast$-homomorphisms from $\Td$ to 
$(1-q)\cdot\mathcal{T}(E)_{| F}\cdot(1-q)$ and $q\cdot\mathcal{T}(E)_{| F}\cdot q$ 
given by 
$s_i\mapsto \ell(\xi_1')\ell(\xi_{1+i}')\,\ell(\xi_1')^*$ and $s_i\mapsto \ell(\xi_{1+i}')\,q$ 
for $i= 1, 2\,$. \\ \indent 
The compactness of the space $X$ enables to end the proof of this first assertion.   

\medskip\noindent 2) 
Let $v_l\in\mathcal{T}(E)_{|\bar{G}_l}$ be the partial isometry 
$v_l:=\ell(\zeta_{l+1})_{|\bar{G}_l}\cdot\ell(\xi(l))^*_{|\bar{G}_l}$. 
There exists by Lemma 2.4 of \cite{BRR08} a K$_1$-trivial unitary $w_l$ 
in the properly infinite unital \cst-algebra $\mathcal{T}(E)_{|\bar{G}_l}$ 
which has the two requested properties (a) and (b). 

\medskip\noindent 3) One constructs inductively the restrictions $\xi_{|G_k}$ in $E_{|G_k}$. 
Set $\xi_{|G_1}:=\zeta_1$ and assume $\xi_{|G_k}$ already constructed. 
As $\ell(\xi_{|G_k})_{|\bar{G}_k}=z_{k+1}^*\cdot\ell(\zeta_{k+1})_{|\bar{G}_k}\,$, 
the only extension $\xi_{|G _{k+1}}\in E_{| G_{k+1}}$ such that 
$(\xi_{|G _{k+1}})_{| G_k}=\xi_{|G_k}$ and 
$(\xi_{|G _{k+1}})_{| F_{k+1}}=\mathbf{\overline{q}}(z_{k+1})^*\cdot(\zeta_{k+1})_{| F_{k+1}}$ 
satisfies $\|(\xi_{|G _{k+1}})_x\|=1$ for all point $x\in G_{k+1}$. 
\qed

\medskip 
\begin{rems} 
a) The nontrivial separable Hilbert $C(B_\infty)$-module $E_{DD}$ 
constructed by J. Dixmier and A. Douady (\cite{DD63}) 
has infinite dimensional fibres and 
every section $\zeta\in E_{DD}$ satisfies $\zeta_x=0$ 
for at least one point $x\in B_\infty$. 
Thus, it cannot satisfy the assumptions for the assertion~3) of Proposition~5.2. 
There are some {$k\in\{1, \ldots, m-1\}$} and 
a unitary  $a_{k+1}\in\U_0\bigl(M_2(\mathcal{T}(E_{DD})_{|F_{k+1}})\bigr)$ such that 
\begin{equation} 
  (a_{k+1})_{|\bar{G}_k}=w_k\oplus 1_{|\bar{G}_k} 
\end{equation} 
and either 
$a_{k+1}\not\in\mathcal{T}(E_{DD})_{|F_{k+1}}\oplus C(F_{k+1})$ or 
 $\alpha(a_{k+1})\not=a_{k+1}\ot 1$. 

 %
\noindent b) If each of the K$_1$-trivial unitaries $w_l$ introduced in assertion 2) 
of Proposition~\ref{prop4.2} satisfies 
$\alpha(w_l)=w_l\ot 1$ and 
$w_l\sim_h 1_{|\bar{G}_l}$ in $\U(\mathcal{T}E)^\alpha{}_{|\bar{G}_l})$, 
then  there exist by \cite[ Lemma~2.1.7]{LLR00} 
$m-1$ unitaries $z_{l+1}\in\mathcal{T}(E)^\alpha{}_{| F_{l+1}}$
such that $(z_{l+1})_{| \bar{G}_l}=w_l$, 
so that there exists a section $\xi\in E$ with $\xi_x\not=0$ for all $x\in X$. 

\noindent c) Let $A$ be a separable unital continuous $C(B_\infty)$-algebra 
with fibres isomorphic to $\Od$ 
such that $K_i(A)\cong C(Y_0, \Gamma_i)$ for $i=0, 1$, 
where $(\Gamma_0, \Gamma_1)$ is a pair of countable abelian torsion groups (\cite[\S 3]{Dad09}). 
Let $\varphi$ be a continuous field of faithful states on $A$. 
Then the $C(B_\infty)$-algebra $A'\subset\L(L^2(A, \varphi) )$ generated by $\pi_\varphi(A)$ and 
the algebra of compact operators $\K(L^2(A, \varphi) )$ is a continuous $C(B_\infty)$-algebra since 
both the ideal $\K(L^2(A, \varphi) )$ and the quotient $A\cong A'/\K(L^2(A, \varphi) )$ are continuous 
(see \eg \cite[Lemma 4.2]{Blan09}). 
All the fibres of $A'$ are isomorphic to the Cuntz extension $\Td$. 
But $A'$ is not a trivial $C(B_\infty)$-algebra since 
$K_0(A')=C(Y_0, \Gamma_0)\oplus\Z$ and $K_1(A')=C(Y_0, \Gamma_1)\,$. 
\end{rems}


\begin{ques} 
The Pimsner-Toeplitz algebra $\mathcal{T}(E_{DD})$ is locally purely infinite 
(\cite[Definition 1.3]{BK04b}) 
since all its simple quotients are isomorphic to the Cuntz algebra $\Oinf$ 
(\cite[Proposition 5.1]{BK04b}). 
But is $\mathcal{T}(E_{DD})$ properly infinite?
\end{ques}

\href{mailto:Etienne.Blanchard@imj-prg.fr}{Etienne.Blanchard@imj-prg.fr} 

\address{IMJ-PRG, \quad UP7D - Campus des Grands Moulins, \quad Case 7012\\ 
${}$\quad F-75205 Paris Cedex 13\\
} 


\begin{thebibliography}{1} 
\bibitem[Blac04]{Blac04} B. Blackadar, 
\textit{Semiprojectivity in simple \cst-algebras}, 
Adv. Stud. Pure Math. \textbf{38} (2004), 1--17. 
\bibitem[Blan97]{Blan97} E. Blanchard, 
\textit{Subtriviality of continuous fields of nuclear C*-algebras}, with an appendix by E. Kirchberg, 
J. Reine Angew. Math. \textbf{489} (1997), 133--149. 
\bibitem[Blan09]{Blan09} E. Blanchard, 
\textit{Amalgamated free products of C*-bundles}, 
Proc. Edinburgh Math. Soc. \textbf{52} (2009), 23--36. 
\bibitem[Blan10]{Blan10} E. Blanchard, 
\textit{K$_1$-injectivity for properly infinite \cst-algebras}, 
Clay Math. Proc. \textbf{11} (2010), 49--54. 
\bibitem[Blan13]{Blan13} E. Blanchard, 
\textit{Continuous fields with fibres $\Oinf$}, 
Math. Scand. \textbf{115} (2014), 189--205.
\bibitem[BK04a]{BK04a} E. Blanchard, E. Kirchberg, 
\textit{Global Glimm halving for C$^*$-bundles}, 
J. Op. Th. \textbf{52} (2004), 385--420. 
\bibitem[BK04b]{BK04b} E. Blanchard, E. Kirchberg, 
\textit{Non-simple purely infinite C*-algebras: the Hausdorff case}, 
J. Funct. Anal. \textbf{207} (2004), 461--513. 
\bibitem[BRR08]{BRR08} E. Blanchard, R. Rohde, M. R\o rdam, 
\textit{Properly infinite $C(X)$-algebras and $K_1$-injectivity}, 
J. Noncommut. Geom. \textbf{2} (2008), 263--282. 
\bibitem[CEI08]{CEI08}K. T. Coward, G. Elliott, C. Ivanescu, 
\textit{The Cuntz semigroup as an invariant for \cst-algebras}, 
J. Reine Angew. Math. \textbf{623} (2008), 161--193. 
\bibitem[Cun77]{Cun77}J. Cuntz, 
\textit{Simple \cst-Algebras Generated by Isometries}, 
Commun. Math. Phys. \textbf{57} (1977), 173--185. 
\bibitem[Cun81]{Cun81}J. Cuntz, 
\textit{K-theory for certain \cst-algebras}, 
Ann. of Math. \textbf{113} (1981), 181Ð-197.
\bibitem[DD63]{DD63} J. Dixmier, A. Douady, 
\textit{Champs continus d'espaces hilbertiens et de C$^*$-alg\`ebres}, 
 Bull. Soc. Math. France \textbf{91} (1963), 227--284. 
\bibitem[Dix69]{Dix69} J. Dixmier, 
\textit{Les \cst-alg\`ebres et leurs repr\'esentations}, 
Gauthiers-Villars Paris (1969).
 \bibitem[D\u{a}d09]{Dad09} M. D\u{a}d\u{a}rlat. 
 \textit{Fiberwise $KK$-equivalence of continuous fields of \cst-algebras}, 
 J. K-Theory \textbf{3} (2009), 205--219. 
\bibitem[DW08]{DW08} M. D\u{a}d\u{a}rlat, W. Winter, 
\textit{ Trivialization of $C(X)$-algebras with strongly self-absorbing fibres}, 
Bull. Soc. Math. France \textbf{136} (2008), 575-Ð606. 
\bibitem[Ell94]{Ell94}G. Elliott, 
 \textit{The classification problem for amenable \cst-algebras}, 
Proc. Internat. Congress of Mathematicians (Zurich, Switzerland, 1994), 
Birkhauser Verlag, Basel (1995), 922--932. 
\bibitem[HRW07]{HRW07}I. Hirshberg, M. R\o rdam, W. Winter, 
\textit{$C_0(X)$-algebras, stability and strongly self-absorbing \cst-algebras}, 
Math. Ann. \textbf{339} (2007), 695--732. 
\bibitem[HR98]{HR98}J. Hjelmborg, M. R\o rdam,
\textit{On stability of \cst-algebras}. 
J. Funct. Anal. \textbf{155} (1998), 153--170. 
\bibitem[Kas88]{Kas88} G.G. Kasparov, 
\textit{Equivariant KK-theory and the Novikov conjecture}, 
Invent. Math. \textbf{91} (1988), 147--201. 
\bibitem[KR00]{KR00}E. Kirchberg, M. R\o rdam, 
\textit{Non-simple purely infinite \cst-algebras}, 
Amer. J. Math. \textbf{122} (2000), 637--666. 
\bibitem[LLR00]{LLR00} F. Larsen, N. J. Laustsen, M. R\o rdam, 
\textit{An Introduction to K-theory for $C^\ast$-algebras}, 
London Mathematical Society Student Texts \textbf{49} (2000) CUP, Cambridge. 
\bibitem[Pim95]{Pim95} M. V. Pimsner, \textit{A class of \cst-algebras generalizing 
both Cuntz-Krieger algebras and crossed products by $\Z$}, 
Free probability theory, 
Fields Inst. Commun. \textbf{12} (1997), 189--212. 
\bibitem[RR11]{RR11}L. Robert, M. R\o rdam, \textit{Divisibility properties for \cst-algebras}, 
Proc. London Math. Soc. \textbf{106} (2013), 1330--1370.
\bibitem[Roh09]{Roh09}R. Rohde, \textit{K$_1$-injectivity of \cst-algebras}, 
Ph.D. thesis at Odense (2009) 
\bibitem[R\o r97]{Ror97} M. R\o rdam, \textit{Stability of \cst-algebras is not a stable property.}, 
Doc. Math. J. \textbf{2} (1997), 375--386. 
\bibitem[R\o r03]{Ror03} M. R\o rdam, 
\textit{A simple \cst-algebra with a finite and an infinite projection}, 
Acta Math. \textbf{191} (2003), 109Ð142. 
\bibitem[R\o r04]{Ror04} M. R\o rdam, \textit{The stable and the real ranks of 
$\mathcal{Z}$-absorbing \cst-algebras}, 
Internat. J. Math. \textbf{15} (2004), 1065--1084. 
\bibitem[TW07]{TW07} A. Toms, W. Winter, \textit{Strongly self-absorbing \cst-algebras}, 
Trans. Amer. Math. Soc. \textbf{359} (2007), 3999Ð-4029. 
\bibitem[Win09]{Win09} W. Winter, 
\textit{Strongly self-absorbing \cst-algebras are $\mathcal{Z}$-stable}, 
J. Noncommut. Geom. \textbf{5} (2011), 253Ð-264. 

\end{thebibliography}
\end{document}